# Shell-infill composite structure design based on a hybrid explicit- implicit topology optimization method

Yilin Guo[1], Chang Liu[1,2*], Xu Guo[1,2*]

[1]*State Key Laboratory of Structural Analysis, Optimization and CAE Software for Industrial Equipment,*
*Department of Engineering Mechanics,*
*International Research Center for Computational Mechanics,*
*Dalian University of Technology, Dalian, 116023, P.R. China*

[2]*Ningbo Institute of Dalian University of Technology, Ningbo, 315016, P.R. China*

**Abstract**

The present paper introduces a hybrid explicit-implicit topology optimization method for shell-infill composite structure design. The proposed approach effectively combines the advantages of the explicit Moving Morphable Component (MMC) method, which describes structural topology only using a set of geometric parameters, and the implicit Solid Isotropic Material with Penalization (SIMP) method, which offers greater design freedom for characterizing the structural features. Compared to the existing methods for shell-infill structure design, the proposed approach can obtain optimized shell-infill structures with complex infill topology without resorting to complex filtering/projection operations. Numerical examples demonstrate the effectiveness of the proposed approach.

**Keywords:** Shell-infill composite structure; Topology optimization; Moving morphable component (MMC); Solid isotropic material with penalization (SIMP).

---

[*] Corresponding authors. E-mail: c.liu@dlut.edu.cn (Chang Liu), guoxu@dlut.edu.cn (Xu Guo)

# 1 Introduction

Porous infill structures offer several advantages over fully solid structures, including a high strength-to-weight ratio [1], efficient energy-absorbing capabilities [2], increased resistance to buckling [3], and improved performance stability under varying loads and material defects [4]. In fact, various lightweight and high-strength shell-infill composite structures can be found in nature, such as animal bones and plant stems. These composite materials typically consist of a rigid outer shell and a porous inner infill, as illustrated in Fig. 1. Furthermore, advancements in modern manufacturing technology, such as additive manufacturing, have made it possible to manufacture shell-infill composite structures (see in Fig. 2) and have facilitated their successful applications in various engineering fields, including automotive, robotics, aerospace, and many others.

In order to make more effective use of shell-infill structures, it is highly desirable to develop efficient design methods. Therefore, in recent years, using topology optimization methods [9–15] to design shell-infill composite structures has become a hot research topic [16–18]. The design and optimization of shell-infill composite structures involves two key issues. The first is determining the morphology of the coated shell, while the second is identifying the topology of the porous infill structure within the coated shell. Various topology optimization methods have been developed to tackle the above challenging issues, and each of them has its unique advantages.

Early studies mainly investigated porous infill structures with uniform infill distributions. In these methods, homogenization methods based on scale-separation and period-distribution assumptions were employed to estimate the effective material properties of the infill structure, thereby reducing the computational cost in the optimization process [19,20]. However, such uniform filling patterns highly restrict the design flexibility and hinder the full potential of using modern manufacturing techniques. Therefore, non-uniform filling structures with spatial varying infill patterns have attracted increasing attention. In general, non-uniform filled structures may offer greater design freedom and therefore achieve more better structural performances. Designing non-uniform filling structures, however, may significantly increases computational efforts. Additionally, ensuring connectivity between different cells in the filling structure during the design process is also a challenge that needs to be addressed.

In order to reduce the computational cost in the optimization of non-uniform porous infilled structures, several approaches have been proposed, including the multi-domain approach [21] and the projection-based process [22]. Furthermore, to ensure the smooth connectivity of non-uniform infill microstructures, the so-called de-homogenization method [23,24], shape metamorphosis technology [25], and other related approaches [26] have been proposed. In contrast to the above methods where the assumption of scale separation has been adopted, there also exist some methods designing of porous-infilled structures directly at the macroscopic scale [4,27,28]. For example, Wu et al. [4] used the SIMP method to perform topology optimization of bone-like porous infill structures within a given region by imposing some macroscopic local material volume constraints. They found that the non-uniform and interconnected infill microstructures aligned with the principal stress directions exhibit higher stiffness than uniform infill. The main advantage of this method is its ability to guarantee the structural connectivity of spatially varying infill structures without introducing additional constraints.

Most of the aforementioned studies have primarily focused on designing porous infill materials. However, it is worth noting that in additive manufacturing, porous infill materials are typically enclosed or protected by a shell. These shells can significantly enhance the load-bearing capacity of composite structures. The design of coated shells is often considered challenging for density-based topology optimization methods. Nevertheless, Clausen et al. [29] successfully addressed this issue by optimizing the topology of coated structures using the classical SIMP (Solid Isotropic Material with Penalization) method. By employing a two-step filtering and projection strategy as well as leveraging spatial gradients at the material interfaces, the uniform-thickness coating shell is identified from the base structure (a crucial step in the optimization process).

Initially, this method was developed for 2D cases and later extended to 3D cases [30]. Subsequently, Luo et al. [31] further simplified this approach and devised a more efficient erosion-based method for identifying the interface between the shell and the base structure in density-based approaches. Additionally, to provide more details on the infill structure, they developed a specific erosion operation to separate the shell and infill components through a two-step filtering and projection process.

In addition to the density-based approach, recent years have also witnessed the application of explicit topology optimization methods, such as Moving Morphable Components (MMC) [15] and Moving Morphable Voids (MMV) methods [32], for the design of shell-infill structures. One of the advantages of explicit optimization methods is their ability to control the topology and geometric features of a structure using only a few design variables. Liu et al. introduced a hybrid MMC-MMV approach for the optimal design of shell-graded-infill structures [33]. In this approach, morphable voids are used to represent the boundary of the shell, while morphable components are employed to depict the distribution of infill material.

The above-mentioned studies have either adopted geometric explicit methods like MMC/MMV or geometric implicit methods like SIMP as the sole topology optimization approach. However, both of these methods have their own advantages and disadvantages. The geometry implicit topology optimization method based on density can achieve more complex topology for infill microstructures but faces some challenges, particularly in simplifying the optimization process for shell structures and generating smooth and geometrically explicit shells. On the other hand, the geometry explicit topology optimization method can provide benefits such as more explicit boundary geometry descriptions and seamless integration with CAD/CAE. However, it cannot generate infill microstructure topologies with fine structural features (which may be very demanded for multi-physics applications) as in the implicit method.

To address the aforementioned issues, this paper introduces a hybrid explicit-implicit topology optimization method for shell-infill structure design. The proposed approach effectively combines the advantages of the explicit Moving Morphable

Component (MMC) method, which describes structural topology using only a set of geometric parameters, and the implicit Solid Isotropic Material with Penalization (SIMP) method, which offers greater design freedom to characterize structural details. Compared to existing methods for shell-infill structure design, the proposed approach can achieve optimized shell-infill structures with complex infill topology without the need for complex filtering/projection operations. Numerical examples demonstrate the effectiveness of the proposed approach.

The remaining sections of the paper are organized as follows: In Section 2, we introduce the hybrid MMC-SIMP method for describing shell-infill structures. Section 3 presents the mathematical formulation of the optimization problem under consideration. A discussion on finite element (FE) discretization, numerical implementation, and sensitivity analysis is provided in Section 4. In Section 5, we validate the effectiveness of the proposed approach through three representative numerical examples. Finally, some concluding remarks are presented.

## 2 Description of shell-infill structures by a hybrid explicit-implicit topology optimization method

This section presents a method to generate shell-infill structures based on explicit-implicit combined topology optimization. The target optimization structure, shown in Fig. 3, consists of a shell with a fixed thickness $t$ and a porous infill structure occupying the domain enclosed by the shell. The shell and infill structure are optimized simultaneously in the proposed algorithm. The shell maintains its thickness $t$ during optimization, while the infill structure forms non-uniform porous structures within the evolving domain defined by the shell. The implicit Solid Isotropic Material with Penalization (SIMP) method is employed to design the porous infill structure, as it can capture the local details of the structure effectively (especially at high background mesh resolution). However, since it lacks explicit geometric information about the structure boundaries, it is not easy to carry out the design of the shell structure, so in this work, we choose the explicit topology optimization method based on the Moving Morphable

Components (MMC) to design of the coated shell structure, and finally, establishes a unified framework for the simultaneous optimization of both parts of the structure.

## 2.1 Generation of infill structures based on SIMP method

The implicit SIMP method is used to generate the porous infill structure. The SIMP method uses element densities as the basic variables of the structure topology description and can capture the complex infill structure finely. The design domain is discretized into a series of finite elements, each with an artificial material density (denoted as $\rho_e^S$ for the $e$th element), as illustrated in Fig. 4. The structure topology can be changed by driving the element density from 0 to 1 through the optimization algorithm. The design domain is firstly discretized into $NE$ finite elements, and then the design variable vector of the SIMP method can be written as $\boldsymbol{D}^{\mathrm{SIMP}} = \left(\rho_1^S, \ldots, \rho_e^S, \ldots, \rho_{NE}^S\right)^\top$. To ensure the generation of clear black-and-white infill structures, an improved SIMP method is adopted [34]. The improved SIMP method can reduce the gray elements in the optimization results by applying additional filtering and projection operations, as follows.

*Filtering.* First, define a collection $\mathcal{M}_e$, which represents the set of elements whose distance to the center of $e$th element (denoted as $\Delta(e, i)$) is less than the prescribed filter radius $r_{\min}$. Next, we define a weight operator $F_{ei}$:

$$F_{ei} = \max\bigl(0, r_{\min} - \Delta(e, i)\bigr), \qquad (2.1)$$

by utilizing the operator $F_{ei}$, the density $\rho_e^S$ is transformed to

$$\tilde{\rho}_e^S = \frac{1}{\sum_{i \in M_e} F_{ei}} \sum_{i \in M_e} F_{ei}\, \rho_i^S. \qquad (2.2)$$

For further information about the filtering operations, please refer to reference [35,36].

*Projection.* In order to ensure the optimization result of the infill structure is either 0 or 1, the intermediate density $\tilde{\rho}_e^S$ can be threshed at the value of 1/2 as follows:

$$\breve{\rho}_e^S(\tilde{\rho}_e^S) = \begin{cases} 1, & \text{if } \tilde{\rho}_e^S > \frac{1}{2} \\ 0, & \text{oterwise} \end{cases} \quad (2.3)$$

For numerical optimization, we relax $\breve{\rho}_e^S$ to a threshold function and approximate this non-differential function by:

$$\breve{\rho}_e^S(\tilde{\rho}_e^S) = \frac{\tanh\left(\frac{\beta}{2}\right) + \tanh\left(\beta\left(\tilde{\rho}_e^S - \frac{1}{2}\right)\right)}{2\tanh\left(\frac{\beta}{2}\right)} \quad (2.4)$$

The parameter $\beta$ is used to control the sharpness of the threshold function. When $\beta$ approaches infinity, the threshold function becomes stricter, resulting in a stringent binary classification. To avoid highly nonlinear situations in the optimization process, $\beta$ is set to 1 in the initial optimization iteration step, and after a certain number of iterations (such as 100 steps), the $\beta$ value is doubled, instead of directly using a larger $\beta$ value. This progressive increase of $\beta$ method helps to maintain a relatively smooth optimization process [37].

Under the SIMP framework, it is convenient to generate non-uniform porous infill structures. In this work, we adopted the method proposed by Wu et al. [4], which introduces a local volume constraint in the optimization formulation. This constraint controls the material usage near each point in the design domain, compelling the material to form a porous infill structure. The mathematical expression for the local volume constraint can be written as:

$$\bar{v}_e = \frac{\sum_{i\in\mathcal{N}_e} \breve{\rho}_i^S A_i}{\sum_{i\in\mathcal{N}_e} A_i} < \bar{V}_{\text{loc}}, e = 1,\ldots,i,\ldots,NE, \quad (2.5)$$

where $\bar{v}_e$ represents the local material volume near the $e$th element, the filtered and projected element density is denoted as $\breve{\rho}_i^S$, symbol $\bar{V}_{\text{loc}}$ denotes the upper bounds of the local volume, symbol $\mathcal{N}_e$ is the collection of all elements within the neighborhood with a radius of $R_e$ centered at the $e$th element and the area of the $i$th finite element can be represented as $A_i$. The definition of $\mathcal{N}_e$ is as follows:

$$\mathcal{N}_e = \{i \mid \|\,\boldsymbol{x}_i - \boldsymbol{x}_e\,\|_2 \leq R_e\}, \quad (2.6)$$

where, $x_i$ and $x_e$ are the coordinates of the center points of the corresponding elements. To avoid numerical issues caused by a large number of constraints introduced by directly using Eq. (2.5) in the optimization formulation, in the numerical implementation, the local volume constraint is further condensed as follows:

$$\max_{1 \leq e \leq NE}(\bar{v}_e) - \bar{V}_{\text{loc}} \leq 0. \tag{2.7}$$

By introducing the local volume constraints as described above into the optimization formulation, it can generate non-uniform porous infill structures as illustrated in Fig. 5.

**2.2 Generation of shell structures based on MMC method**

The pixel-based topology description of the SIMP method is convenient for generating complex porous infill structures. However, it lacks explicit geometric information about the structural boundary, which makes it challenging to construct coated shell structures. As shown by existing work, generating coated shell structures under the SIMP framework is tedious. In contrast, the Moving Morphable Components (MMC) explicit topology optimization method, developed in recent years, uses a group of components with explicit geometric parameters as the building blocks of structure topology description and performs the topology optimization by moving, deforming, overlapping, and merging the components. The MMC method naturally has explicit geometric parameters of the structural boundary and can generate coated shell structures conveniently through intuitive geometric ways.

In the MMC method, the topology description function (TDF) is employed to characterize the domain occupied by the components. The TDF of whole structure is defined as follows:

$$\begin{cases} \varphi^s(x) > 0, & \text{if } x \in \Omega^s, \\ \varphi^s(x) = 0, & \text{if } x \in \partial\Omega^s, \\ \varphi^s(x) < 0, & \text{if } x \in D\setminus(\Omega^s \cup \partial\Omega^s), \end{cases} \tag{2.8}$$

where $D$ represents the prescribed design domain and $\Omega^s \in D$ represents the region occupied by solid materials.

As shown in Fig. 6, the region $\Omega^s$ is composed of $NC$ solid components and the TDF $\varphi^s(x)$ can be obtained by the following formula:

$$\varphi^s(x) = \max\left(\phi_1^C(x), \ldots \phi_i^C(x), \ldots, \phi_{NC}^C(x)\right), \qquad (2.9)$$

where $\phi_i^C(x)$ denotes the TDF of the $i$th component, which is defined as:

$$\phi_i^C(x) = 1 - \left(\frac{x'}{a_i}\right)^p - \left(\frac{y'}{b_i(x')}\right)^p, \qquad (2.10)$$

with

$$\begin{Bmatrix} x' \\ y' \end{Bmatrix} = \begin{pmatrix} \cos\theta_i & \sin\theta_i \\ -\sin\theta_i & \cos\theta_i \end{pmatrix} \begin{Bmatrix} x - x_{0i} \\ y - y_{0i} \end{Bmatrix}, \qquad (2.11)$$

and $p$ is a relatively large positive even number (e.g. $p = 6$).

In Eq. (2.10) and Eq. (2.11), the symbols $(x_{0i}, y_{0i})$, $a_i$, $b_i(x')$ and $\theta_i$ denote the coordinates of the center, the half-length, the variable half-width and the inclined angle (measured from the horizontal in axis anti-clockwise direction) of the $i$-th component, respectively, as shown in Fig. 7.

It should also be noted that the variation of the width $b_i(x')$ is measured with respect to a local coordinate system and $b_i(x')$ can take different forms [38]. In the present work, it is chosen as:

$$b_i(x') = \frac{t_i^1 + t_i^2}{2} + \frac{t_i^2 - t_i^1}{2a_i} x', \qquad (2.12)$$

where $t_i^1$ and $t_i^2$ are parameters that used to describe the thickness of the $i$-th component, as shown in Fig. 7.

The design variables for the MMC method are represented as $\boldsymbol{D}^{\text{MMC}} = ((\boldsymbol{D}_1)^\top, \ldots, (\boldsymbol{D}_i)^\top, \ldots, (\boldsymbol{D}_{NC})^\top)^\top$, with $\boldsymbol{D}_i = (x_{0i}, y_{0i}, a_i, t_i^1, t_i^2, \theta_i)^\top$ represents the design parameters of the $i$th component.

In the numerical implementation, the TDF can be linked to the analytical model through a ersatz material model [39,40]. In this work, to maintain consistency with the generation of porous infill structures based on the SIMP method, the design domain is discretized using the same elements as in Section 2.1, and the density of the $e$th element can be expressed through the surrogate material model as:

$$\rho_e^{\text{MMC}} = \frac{\sum_{i=1}^{4} H_\epsilon((\varphi^s)_i^e)}{4}, e = 1,2,\dots,NE, \quad (2.13)$$

where $H_\epsilon(x)$ is the regularized Heaviside function, the symbol $(\varphi^s)_i^e$ is the TDF value of the $i$th node of the $e$th element. The regularized Heaviside function used in this paper is given by the following expression:

$$H_\epsilon(x) = \begin{cases} 1, & \text{if } x > \epsilon, \\ \frac{3}{4}\left(\frac{x}{\epsilon} - \frac{x^3}{3\epsilon^3}\right) + \frac{1}{2}, & \text{if } -\epsilon \leq x \leq \epsilon, \\ 0, & \text{otherwise}, \end{cases} \quad (2.14)$$

where $\epsilon$ a is a small positive number used to control the length of the transition zone between the solid and void parts.

In the MMC framework, the coated shell structure can be constructed quite simply. Take the structure formed by two components shown in Fig. 8 (a) as an example, assuming that the solid area occupied by the components is $\Omega^{\text{ext}}$. By reducing the length and width of the two components by $2t$ (since the length and width parameters of the components are design variables, this operation can be realized by modifying the corresponding parameters of the components), a new set of components can be obtained, and the solid area occupied by the new components is $\Omega^{\text{int}}$, as shown in Fig. 8(b). Next, a simple Boolean operation on the above two structures will result in a coated shell structure $\Omega^{\text{shell}} = \Omega^{\text{ext}} \backslash \Omega^{\text{int}}$, as shown in Fig. 8(c).

Compared with the tedious construction process of coated shell structures under the SIMP framework, coated shell structures can be conveniently generated in a geometrically intuitive way under the MMC explicit topology optimization framework, and for the equal-thickness coated shell, only one coating thickness parameter needs to be introduced.

**2.3 Generation of shell-infill structures based on MMC-SIMP combined method**

Although the geometric description of porous infill structure and coated shell structure can be realized based on the MMC method and SIMP method, respectively, there is still a problem to be solved, that is, how to combine the two methods to achieve a holistic geometric description of a shell-infill structure. For this purpose, we use the design variables related to MMC and SIMP to combine a new artificial material density vector $\boldsymbol{\rho}$, which should satisfy: $\rho_e = \breve{\rho}_e^S$, when the $e$th element is in $\Omega^{\text{int}}$; $\rho_e = 1$, when the $e$th element is in $\Omega^{\text{shell}}$; $\rho_e = 0$, otherwise.

To achieve this, we first discretize the design domain into $NE$ finite elements. Using the design variables $\boldsymbol{D}^{\text{MMC}}$ related to the MMC method introduced in Section 2.2, and following the method described in Section 2.2, we obtain two sets of components with characteristic size differences $t$, occupying regions $\Omega^{\text{ext}}$ and $\Omega^{\text{int}}$, respectively. From $\boldsymbol{D}^{\text{MMC}}$, we can generate artificial material density vectors $\boldsymbol{\rho}^1 \in \mathrm{R}^{NE}$ and $\boldsymbol{\rho}^2 \in \mathrm{R}^{NE}$ corresponding to each finite element (as shown in Fig. 9 (a) and Fig. 9 (b). The elements in the vector $\boldsymbol{I} - \boldsymbol{\rho}^2$（$\boldsymbol{I} = (1,\dots,1)^\top \in \mathrm{R}^{NE}$）that are not equal to zero correspond to the dual structure of $\Omega^{\text{int}}$: $\overline{\Omega^{\text{int}}} = \mathrm{D}\backslash\Omega^{\text{int}}$ (where $D$ represents the entire design domain). Fig. 9 (d) shows that the white region outside $\overline{\Omega^{\text{int}}}$ is identical to the region that should be filled with porous structure (red part in Fig. 9 (c)).

Similarly, utilizing the design variables $\boldsymbol{D}^{\text{SIMP}}$ from the SIMP method and applying filtering and projection operations, we can obtain another element density vector $\breve{\boldsymbol{\rho}}^S$. The region corresponding to non-zero elements in $\breve{\boldsymbol{\rho}}^S$ can be denoted as $\Omega^{\text{SIMP}}$. Using $\Omega^{\text{SIMP}}$ and $\overline{\Omega^{\text{int}}}$, we can obtain a region $\Omega^1 = \Omega^{\text{SIMP}} \cup \overline{\Omega^{\text{int}}}$.

If we define

$$\boldsymbol{\rho}' = \max(\breve{\boldsymbol{\rho}}^S, \boldsymbol{I} - \boldsymbol{\rho}^2), \tag{2.15}$$

then the region $\Omega^1$ corresponds to the non-zero elements of the vector $\boldsymbol{\rho}'$, as shown in Fig. 9 (e).

If we define

$$\boldsymbol{\rho} = \min(\boldsymbol{\rho}', \boldsymbol{\rho}^1) \ . \tag{2.16}$$

then the region $\Omega^{S-I} = \Omega^1 \cap \Omega^{ext}$ corresponds to the non-zero elements of the vector $\boldsymbol{\rho}'$ (as shown in Fig. 9 (f)) and it is the spatial area occupied by the shell-infill structure that we want to obtain.

In summary, we have

$$\boldsymbol{\rho} = \min(\max(\tilde{\boldsymbol{\rho}}^S, \boldsymbol{I} - \boldsymbol{\rho}^2), \boldsymbol{\rho}^1). \tag{2.17}$$

Compared to implicit methods like density-based approaches that require multiple filtering/projection operations and complex density transformations to construct the topological description for shell-infill structures, this approach fully leverages the flexibility and ease of geometric configuration description offered by the MMC method. It only requires a single-density transformation to obtain the geometric description of the shell-infill structure, greatly simplifying the analytical representation of structural configurations. Additionally, since the geometric details of the internal infill structure are described by the SIMP method, this approach, when compared to using MMC/MMV alone, provides more degrees of freedom for the topological design of the infill structure.

## 3 Problem formulation

In the present study, the shell-infill structures are designed to minimize the structural compliance under the local and total volume constraints of available solid material. Under this circumstance, the corresponding problem formulation can be expressed as follows:

$$\text{Find } \boldsymbol{D} = ((\boldsymbol{D}^{\text{MMC}})^\top, (\boldsymbol{D}^{\text{SIMP}})^\top)^\top$$

$$\text{Minimize } \tilde{C} = C + \gamma C^{\text{M}} = \boldsymbol{u}^\top \boldsymbol{K} \boldsymbol{u} + \gamma (\boldsymbol{u}^C)^\top \boldsymbol{K}^C \boldsymbol{u}^C$$

S. t.

$$\boldsymbol{K}\boldsymbol{u} = \boldsymbol{f},$$

$$\boldsymbol{K}^C \boldsymbol{u}^C = \boldsymbol{f},$$

$$g_1 = \sum_{e=1}^{NE} \rho_e^1 \frac{A_e}{|D|} - \bar{V}_c \leq 0,$$

$$g_2 = \max_{1 \leq e \leq NE}(\bar{v}_e) - \bar{V}_{\text{loc}} \leq 0,$$

$$\boldsymbol{D}^{\text{MMC}} \in \mathcal{D}_{\text{MMC}}, \quad \boldsymbol{D}^{\text{SIMP}} \in \mathcal{D}_{\text{SIMP}}, \tag{3.1}$$

where $\boldsymbol{u}$ is the displacement vector of the structure, $\boldsymbol{f}$ is the vector of the external load and $\boldsymbol{K}$ is the global stiffness matrix of the structure (see the discussions in the next section). In Eq. (3.1), $C = \boldsymbol{u}^\top \boldsymbol{K} \boldsymbol{u}$ represents the structural compliance while $C^M = (\boldsymbol{u}^C)^\top \boldsymbol{K}^C \boldsymbol{u}^C$ represents the compliance of the structure formed solely by lager-size components (as shown in Fig. 9(a)) under the external load with $\boldsymbol{K}^C$ and $\boldsymbol{u}^C$ denoting the corresponding global stiffness matrix and the displacement vector, respectively. The integration of $C^M$ into the objective function is to promote the fusion of the morphable components (see the discussions in Section 5 for the necessity of this operation). In addition, $g_1$ is the volume constraint function associated with the morphable components and $g_2$ represents the aggregated local volume constraint function (see the discussion in Section 2), respectively. In Eq. (3.1), $\gamma$ is a parameter to manifest the effect of the penalty term $C^M$, $A_e$ is the area of the $e$th element and $|D|$ is the total area of the design domain. The symbols $\bar{V}_c$ and $\bar{V}_{\text{loc}}$ are the upper bounds of two volume constraints, respectively. Furthermore, $\mathcal{D}_{\text{MMC}}$ and $\mathcal{D}_{\text{SIMP}}$ denote the feasible sets of the MMC-type and SIMP-type design variables, respectively.

## 4 Numerical implementation aspects and sensitivity analysis

### 4.1 Finite element analysis

In this study, the design domain is discretized into uniform four-node quadrilateral plane stress elements (Q4 elements). In this context, the stiffness matrix of the $e$th finite element can be calculated as follow:

$$\boldsymbol{k}_e(\rho_e) = \boldsymbol{k}_{\min} + (\rho_e)^n (\boldsymbol{k}^S - \boldsymbol{k}_{\min}), \tag{4.1}$$

where n is a penalty factor (we take $n = 3$ in this paper), $\boldsymbol{k}^S$ is the element stiffness matrix of Q4 element when $\rho_e = 1$ and $\boldsymbol{k}_{\min} = \boldsymbol{k}(\rho_{\min})$ with $\rho_{\min}$ is small lower bound of the material density introduced for avoiding the singularity of the global stiffness matrix (we take $\rho_{\min} = 10^{-3}$). Global stiffness matrix $\boldsymbol{K}$ is assembled from the stiffness matrix $\boldsymbol{k}_e(e = 1,2,...NE)$ of each finite element.

**4.2 Sensitivity analysis**

Based on the ersatz material model adopted in finite element analysis, we can obtain the sensitivity of the objective function and constraint functions with respect to an arbitrary design variable $d$ analytically. To this end, we first employ the K-S function [41] to approximate the max and min operators involved in the corresponding calculations

$$\chi = \frac{\ln(\sum_{i=1}^{m} \exp(L\chi^i))}{L}. \tag{4.2}$$

Actually, when $L$ is a relatively large negative integer (e.g., $L_1 = -100$), it yields that $\chi \approx \min(\chi^1,...\chi^i,...,\chi^m)$, and when $L$ is a relatively large positive integer (e.g., $L_2 = 200$), we have $\chi \approx \max(\chi^1,...,\chi^i,...\chi^m)$.

The sensitivity of the structural compliance with respect to a typical design variable $d$ can be expressed as:

$$\frac{\partial C}{\partial d} = \boldsymbol{u}^\top \frac{\partial \boldsymbol{K}}{\partial d} \boldsymbol{u} = \sum_{e=1}^{NE} -\boldsymbol{u}_e^\top \frac{\partial \boldsymbol{k}_e}{\partial d} \boldsymbol{u}_e, \tag{4.3}$$

where $\boldsymbol{k}_e$ and $\boldsymbol{u}_e$ are the element stiffness matrix and displacement vector of the $e$-th element, respectively.

As for the calculation of $\partial \boldsymbol{k}_e/\partial d$, the following two cases should be considered.

(1) If $d \in \boldsymbol{D}^{\text{MMC}}$,

$$\frac{\partial \boldsymbol{k}_e}{\partial d} = \frac{\partial \boldsymbol{k}_e}{\partial \rho_e}\left(\frac{\partial \rho_e}{\partial \rho_e'}\frac{\partial \rho_e'}{\partial \rho_e^2}\frac{\partial \rho_e^2}{\partial d} + \frac{\partial \rho_e}{\partial \rho_e'}\frac{\partial \rho_e'}{\partial \rho_e^1}\frac{\partial \rho_e^1}{\partial d}\right) = \frac{n}{4}\sum_{e=1}^{NE}\sum_{i=1}^{4} G_1(\rho_e)^{n-1}\,\boldsymbol{k}^S, \tag{4.4}$$

where

$$G_1 = \frac{A+B}{\exp(L_1\rho_e') + \exp(L_1\rho_e^1)}, \tag{4.5}$$

$$A = -\exp(L_1\rho'_e)\frac{\exp(L_2(1-\rho_e^2))}{\exp(L_2\rho_e^S) + \exp(L_2(1-\rho_e^2))}\frac{\partial H((\varphi^S)_i^e)}{\partial d}, \quad (4.6)$$

$$B = \exp(L_1\rho_e^1)\frac{\partial H((\varphi'^S)_i^e)}{\partial d}, \quad (4.7)$$

where the values of $\partial H((\varphi^S)_i^e)/\partial d$ and $\partial H((\varphi'^S)_i^e)/\partial d$ can be obtained by the method provided in [40]. Then we have:

$$\frac{\partial C}{\partial d} = -\frac{n}{4}\sum_{e=1}^{NE}\sum_{i=1}^{4} G_1(\rho_e)^{n-1}\boldsymbol{u}_e^\top \boldsymbol{k}^S \boldsymbol{u}_e. \quad (4.8)$$

(2) If $d \in \boldsymbol{D}^{\text{SIMP}}$,

$$\frac{\partial \boldsymbol{k}_e}{\partial d} = \frac{\partial \boldsymbol{k}_e}{\partial \rho_e}\frac{\partial \rho_e}{\partial \rho'_e}\frac{\partial \rho'_e}{\partial \check{\rho}_e^S}\frac{\partial \check{\rho}_e^S}{\partial \tilde{\rho}_e^S}\frac{\partial \tilde{\rho}_e^S}{\partial \rho_e^S}\frac{\partial \rho_e^S}{\partial d} = n\sum_{e=1}^{NE}(\rho_e)^{n-1} G_2 \frac{\partial \check{\rho}_e^S}{\partial \tilde{\rho}_e^S}\frac{\partial \tilde{\rho}_e^S}{\partial \rho_e^S}\boldsymbol{k}^S, \quad (4.9)$$

where

$$G_2 = \frac{\exp(L_1\rho'_e)}{\exp(L_1\rho'_e) + \exp(L_1\rho_e^1)}\frac{\exp(L_2\rho_e^S)}{\exp(L_2\rho_e^S) + \exp(L_2(1-\rho_e^2))}, \quad (4.10)$$

$$\frac{\partial \check{\rho}_e^S}{\partial \tilde{\rho}_e^S} = \frac{\beta(1-\tanh^2(\beta(\tilde{\rho}_e - \frac{1}{2})))}{2\tanh\left(\frac{\beta}{2}\right)}, \quad (4.11)$$

$$\frac{\partial \tilde{\rho}_e^S}{\partial \rho_e^S} = \sum_{e\in M_j}\frac{1}{\sum_{i\in M_e} F_{ei}}F_{je}. \quad (4.12)$$

Then we have

$$\frac{\partial C}{\partial d} = -n\sum_{e=1}^{NE}\sum_{i=1}^{4}(\rho_e)^{n-1} G_2 \frac{\partial \check{\rho}_e^S}{\partial \tilde{\rho}_e^S}\frac{\partial \tilde{\rho}_e^S}{\partial \rho_e^S}\boldsymbol{u}_e^\top \boldsymbol{k}^S\boldsymbol{u}_e. \quad (4.13)$$

and the value of $\partial C^M/\partial d$ can also be obtained in a similar way.

The derivation of sensitivities of the constraint functions with respect to design variables is trivial and will not be discussed here.

## 5 Numerical examples

In this section, three plane stress problems are investigated to demonstrate the effectiveness of the proposed method. To simplify the analysis, it is assumed that all relevant quantities in the studied problems are dimensionless, and the thickness of all design domains is set to a unit value. The Young's modulus and Poisson's ratio of the isotropic solid material are chosen as $E_s = 1$ and $v_s = 0.3$, respectively. Additionally, the design variables are updated using the MMA algorithm [42], and the optimization process is terminated if either of the following conditions is met: (1) the maximum relative change in each design variable between two consecutive iterations falls below a specified threshold (i.e., 0.05%), or (2) the number of iteration steps reaches a predefined threshold (i.e., 2000). In all examples, uniform square-shaped plane stress elements are used for finite element analysis. To avoid trapping into a local optimum prematurely, we choose not to refrain from directly apply a larger projection sharpness parameter. Instead, we initiate the optimization process with a lower sharpness value, such as $\beta = 1$. As the local volume constraint approaches a constant value during the optimization, we double the sharpness parameter value to maintain the stability of the optimization process and the accuracy of finite element analysis.

**5.1 The L-shaped beam**

The dimensions of the L-beam and the design domain considered in this example are shown in Fig. 10. The design domain is discretized with square elements at a resolution of 300×300 to better capture the geometric details of the infill structure. For computational convenience, a 1×1 non-designable domain is placed in the upper right corner of the 2×2 design domain, resulting in an L-shaped design domain. The optimization parameters for this example are selected as follows: the density design variable filter radius is set to $r_{\min} = 2a$ (where $a$ represents the size of the finite element), the influence radius for the local volume constraint is $R_e = 8a$, the control parameter for the shell thickness is fixed as $t = 3a$, the upper bounds for the component volume constraint and local volume are $\bar{V}_c = 0.5$ and $\bar{V}_{\text{loc}} = 0.6$, respectively.

The optimization results for the L-beam in this example are shown in Fig. 11. The optimization results shown in Fig. 11(a) were conducted according to the optimization formulation in Eq. (3.1) without the introduction of a regularization term to promote component fusion (i.e., γ=0). It can be observed that multiple components with solid shells appear in the optimization results. However, these incompletely fused components, as shown in Fig. 12(a), are not desirable since they would prevent the formation of a large area of continuous infill structure, as depicted in Fig. 12(b). To address this issue, we introduced a regularization term to promote component fusion into the objective function (setting γ≠0 in Eq. (3.1)). Fig. 11(b) plots the optimization results with γ=0.8 (this value will be used for all problems considered in the following). From the figure, it can be seen that after incorporating the regularization term into the objective function, the optimization algorithm successfully fused the components together effectively. In the obtained optimized configuration, there exists a large area of continuous infill structure.

**5.2 A MBB example**

In this example, we examine the well-known MBB problem. The design domain, external loads, and boundary conditions for this problem are illustrated schematically in Fig. 13. Due to the symmetry of the considered problem, only half of the design domain is discretized using 1200×600 square Q4 elements. The density design variable filter radius in this example is set to $r_{\min} = 2a$, and the control parameter for shell thickness is chosen as $t = 8a$. In this example, we mainly focus on testing the algorithm's performance under different component volume constraints and local volume constraints.

Fig. 14 presents the optimization results under three different component volume constraint settings ($\bar{V}_c$ are set as 0.5, 0.6, and 0.7, respectively), while in the local volume constraints $\bar{V}_{\mathrm{loc}} = 0.6$. From Fig. 14, it can be observed that the optimized structure exhibits two different scales of structural details. At a larger scale, the overall topological configuration of the structure is defined by the moving morphable components. At a smaller scale, the infill structure takes on different forms based on the characteristics of stress distribution in different regions. In regions where compressive/tensile uniaxial stress state dominate, the material distribution of the infill structure is generally aligned with the principal stress direction. In areas where multiple components intersect, due to the presence of complex stress states, the infill structure exhibits corresponding cross-shaped material distribution patterns. These results are in accordance with mechanical intuition and are also consistent with findings reported in

[43]. Furthermore, it is observed that as the upper limit of the component volume increases, the compliance of the optimized structure monotonically decreases. This is attributed to the monotonic dependence of structural compliance on material usage.

Fig. 15 illustrates the correlation between the porous infill structure and the background region composed of components. It can be observed that the infill structure is accurately distributed within the region defined by the moving morphable components. Fig. 16 presents some intermediate results of the optimization process. It is observed that in the initial optimization stage, the components establish the primary force transmission pathways and form a distinct coated shell structure profile. Subsequently, through the filtering and projection process, the configuration of the infill structure gradually becomes clearer. Fig. 17 displays the impact of different local volume influence radii on the optimization results, with all other parameters fixed ($\bar{V}_c = 0.6$). It is found that as the influence radius decreases, the porous infill structure becomes denser. Therefore, the sparsity of the porous filling structure can be controlled by adjusting the value of the influence radius. Fig. 18 provides the iterative curve of the optimization process, which is relatively stable, with slight fluctuations possibly attributed to the variations of the values of projection parameter $\beta$.

**5.3 A multiple loads beam example**

Multiple load cases are common in engineering applications, so in this example, a structure similar to a bridge that withstands multiple load conditions was considered to validate the effectiveness of the proposed approach. In this case, the horizontal and vertical displacements at the lower-left and lower-right points of the rectangular design domain are both set to zero, and the design domain is subjected to five vertical loads as shown in Fig. 19. In this example, we take the average of the structural compliance under five different loads as the objective function. Considering the symmetry of the problem, we only optimize half of the design domain. The half of the design domain involved in the computation is discretized with 300×300 square elements. The filtering radius for density design variables is set to $r_{\min} = 2a$, the influence radius for local volume constraints is set to $R_e = 8a$, the control parameter for shell thickness is fixed at $t = 8$, the upper bounds for the component volume constraint and local volume are $\bar{V}_\text{c} = 0.6$ and $\bar{V}_\text{loc} = 0.6$, respectively.

Fig. 20 displays the optimized shell-filling structure obtained through the optimization algorithm, with a corresponding objective function value of 83.98. It can be observed that the method proposed in this paper still exhibits outstanding applicability under multiple loading conditions. In fact, the presence of multiple loading conditions has led to the emergence of more complex filling patterns in the optimized configuration (as shown in Fig. 21). Specifically, for this problem, we have achieved a shell with uniform thickness, which is based on a single connected region formed by moving morphable components. This highlights the effectiveness of the proposed approach in determining the outer profile of the shell-infill structure. Furthermore, as shown in Figure 21, the infill structure is distributed in nearly every part of the base region and exhibits good compatibility with the outer shell.

## 6 Conclusion

In this paper, we introduce a new framework that combines the advantages of explicit and implicit methods. Within this framework, we have developed a topology optimization method that can simultaneously optimize the shell and infill structures. In this approach, the explicit topology optimization method (MMC) is employed to optimize the shell with a specified thickness, while the implicit topology optimization method is used to optimize the internal infill structure with a given local volume fraction constraint. An objective function with a correction term is utilized by us to facilitate the integration of components, thereby obtaining a shell-infill structure with explicit boundary descriptions and intricate details of the infill structure.

We validated the effectiveness of the proposed method through numerical examples. The optimized shell-infill structure exhibits interesting structural patterns. For instance, regions dominated by uniaxial stress display unidirectional infill distribution, while connection areas of multiple components exhibit cross-shaped infill distribution. These findings are in accordance with our understanding of optimal structures and are consistent with mechanical principles. Furthermore, the optimization process typically begins by forming a complete envelope of the shell-infill structure (whose morphology is continuously adjusted during the subsequent optimization process), and then the internal infill structure gradually emerges in the later stages of optimization. Additionally, extending the proposed method from two dimensions to three dimensions is relatively straightforward because both MMC and SIMP methods can be readily implemented in three-dimensional scenarios.

Although the proposed exhibits promising performances, there is still room for improvement in terms of convergence speed, particularly when optimizing the geometric details of the infill structure. Future works will consider more mechanical performances (e.g., buckling behavior, stress/fatigue constraint), manufacturability requirements and employing additive manufacturing equipment for mechanical experimental analysis to validate the findings and conclusions. Furthermore, since most problems in engineering practice inevitably involve various uncertainties, including manufacturing tolerances, variations in loads, and material property heterogeneity, in the future, we plan to explore the topology optimization design of shell-fill structures considering uncertainty based on the method proposed in this paper, combined with effective approaches suggested by Meng et al. [44-46]. Extending the proposed method from two-dimension to three-dimension case (which is much more computationally expansive) is also on the agenda with the aim of expanding its potential applications.

**Acknowledge**

This work is supported by the National Key Research and Development Plan (2022YFB3303000), the National Natural Science Foundation (11821202, 12002077), the Liaoning Revitalization Talents Program (XLYC2001003), the Fundamental Research Funds for the Central Universities (DUT21RC(3)076), and 111 Project (B14013).

**Data Availability Statement**

The data are available from the corresponding author on reasonable request.

# References

[1] Gibson LJ. Cellular Solids. MRS Bull 2003;28:270–4. https://doi.org/10.1557/mrs2003.79.

[2] Ha CS, Lakes RS, Plesha ME. Design, fabrication, and analysis of lattice exhibiting energy absorption via snap-through behavior. Mater Des 2018;141:426–37. https://doi.org/10.1016/j.matdes.2017.12.050.

[3] Clausen A, Aage N, Sigmund O. Exploiting Additive Manufacturing Infill in Topology Optimization for Improved Buckling Load. Engineering 2016;2:250–7. https://doi.org/10.1016/J.ENG.2016.02.006.

[4] Wu J, Aage N, Westermann R, Sigmund O. Infill Optimization for Additive Manufacturing—Approaching Bone-Like Porous Structures. IEEE Trans Vis Comput Graph 2018;24:1127–40.
https://doi.org/10.1109/TVCG.2017.2655523.

[5] Al-Rukaibawi LS, Omairey SL, Károlyi G. A numerical anatomy-based modelling of bamboo microstructure. Constr Build Mater 2021;308:125036. https://doi.org/10.1016/j.conbuildmat.2021.125036.

[6] Sullivan TN, Wang B, Espinosa HD, Meyers MA. Extreme lightweight structures: avian feathers and bones. Mater Today 2017;20:377–91. https://doi.org/10.1016/j.mattod.2017.02.004.

[7] Sarén M-P, Serimaa R, Andersson S, Paakkari T, Saranpää P, Pesonen E. Structural Variation of Tracheids in Norway Spruce (Picea abies [L.] Karst.). J Struct Biol 2001;136:101–9.

[8] Li H, Li H, Gao L, Zheng Y, Li J, Li P. Topology optimization of multi-phase shell-infill composite structure for additive manufacturing. Eng Comput 2023. https://doi.org/10.1007/s00366-023-01837-4.

[9] Bendsøe MP, Kikuchi N. Generating optimal topologies in structural design using a homogenization method. Comput Methods Appl Mech Eng 1988;71:197–224. https://doi.org/10.1016/0045-7825(88)90086-2.

[10] Bendsøe MP. Optimal shape design as a material distribution problem. Struct Optim 1989;1:193–202. https://doi.org/10.1007/BF01650949.

[11] Wang MY, Wang X, Guo D. A level set method for structural topology optimization. Comput Methods Appl Mech Eng 2003;192:227–46. https://doi.org/10.1016/S0045-7825(02)00559-5.

[12] Eschenauer HA, Olhoff N. Topology optimization of continuum structures: A review. Appl Mech Rev 2001;54:331–90. https://doi.org/10.1115/1.1388075.

[13] Xie YM, Steven GP. A simple evolutionary procedure for structural optimization. Comput Struct 1993;49:885–96. https://doi.org/10.1016/0045-7949(93)90035-C.

[14] Guo X, Cheng G-D. Recent development in structural design and optimization. Acta Mech Sin 2010;26:807–23. https://doi.org/10.1007/s10409-010-0395-7.

[15] Guo X, Zhang W, Zhong W. Doing Topology Optimization Explicitly and Geometrically—A New Moving Morphable Components Based Framework. J Appl Mech 2014;81:081009. https://doi.org/10.1115/1.4027609.

[16] Yu H, Huang J, Zou B, Shao W, Liu J. Stress-constrained shell-lattice infill structural optimisation for additive manufacturing. Virtual Phys Prototyp

2020;15:35–48. https://doi.org/10.1080/17452759.2019.1647488.

[17] Das S, Sutradhar A. Multi-physics topology optimization of functionally graded controllable porous structures: Application to heat dissipating problems. Mater Des 2020;193:108775. https://doi.org/10.1016/j.matdes.2020.108775.

[18] Xu S, Liu J, Huang J, Zou B, Ma Y. Multi-scale topology optimization with shell and interface layers for additive manufacturing. Addit Manuf 2021;37:101698. https://doi.org/10.1016/j.addma.2020.101698.

[19] Wang F, Sigmund O. Numerical investigation of stiffness and buckling response of simple and optimized infill structures. Struct Multidiscip Optim 2020;61:2629–39. https://doi.org/10.1007/s00158-020-02525-3.

[20] Wang Y, Wang MY, Chen F. Structure-material integrated design by level sets. Struct Multidiscip Optim 2016;54:1145–56. https://doi.org/10.1007/s00158-016-1430-5.

[21] Liu C, Du Z, Zhang W, Zhang X, Mei Y, Guo X. Design of optimized architected structures with exact size and connectivity via an enhanced multidomain topology optimization strategy. Comput Mech 2021;67:743–62. https://doi.org/10.1007/s00466-020-01961-8.

[22] Groen JP, Wu J, Sigmund O. Homogenization-based stiffness optimization and projection of 2D coated structures with orthotropic infill. Comput Methods Appl Mech Eng 2019;349:722–42. https://doi.org/10.1016/j.cma.2019.02.031.

[23] Groen JP, Sigmund O. Homogenization-based topology optimization for high-resolution manufacturable microstructures. Int J Numer Methods Eng 2018;113:1148–63. https://doi.org/10.1002/nme.5575.

[24] Kim D, Lee J, Nomura T, Dede EM, Yoo J, Min S. Topology optimization of functionally graded anisotropic composite structures using homogenization design method. Comput Methods Appl Mech Eng 2020;369:113220. https://doi.org/10.1016/j.cma.2020.113220.

[25] Wang Y, Zhang L, Daynes S, Zhang H, Feih S, Wang MY. Design of graded lattice structure with optimized mesostructures for additive manufacturing. Mater Des 2018;142:114–23. https://doi.org/10.1016/j.matdes.2018.01.011.

[26] Du Z, Zhou X-Y, Picelli R, Kim HA. Connecting Microstructures for Multiscale Topology Optimization With Connectivity Index Constraints. J Mech Des 2018;140. https://doi.org/10.1115/1.4041176.

[27] Alexandersen J, Lazarov BS. Topology optimisation of manufacturable microstructural details without length scale separation using a spectral coarse basis preconditioner. Comput Methods Appl Mech Eng 2015;290:156–82. https://doi.org/10.1016/j.cma.2015.02.028.

[28] Dou S. A projection approach for topology optimization of porous structures through implicit local volume control. Struct Multidiscip Optim 2020;62:835–50. https://doi.org/10.1007/s00158-020-02539-x.

[29] Clausen A, Aage N, Sigmund O. Topology optimization of coated structures and material interface problems. Comput Methods Appl Mech Eng 2015;290:524–41. https://doi.org/10.1016/j.cma.2015.02.011.

[30] Clausen A, Andreassen E, Sigmund O. Topology optimization of 3D shell

structures with porous infill. Acta Mech Sin 2017;33:778–91. https://doi.org/10.1007/s10409-017-0679-2.

[31] Luo Y, Li Q, Liu S. Topology optimization of shell–infill structures using an erosion-based interface identification method. Comput Methods Appl Mech Eng 2019;355:94–112. https://doi.org/10.1016/j.cma.2019.05.017.

[32] Zhang W, Chen J, Zhu X, Zhou J, Xue D, Lei X, et al. Explicit three dimensional topology optimization via Moving Morphable Void (MMV) approach. Comput Methods Appl Mech Eng 2017;322:590–614. https://doi.org/10.1016/j.cma.2017.05.002.

[33] Liu C, Du Z, Zhu Y, Zhang W, Zhang X, Guo X. Optimal design of shell-graded-infill structures by a hybrid MMC-MMV approach. Comput Methods Appl Mech Eng 2020;369:113187. https://doi.org/10.1016/j.cma.2020.113187.

[34] Andreassen E, Clausen A, Schevenels M, Lazarov BS, Sigmund O. Efficient topology optimization in MATLAB using 88 lines of code. Struct Multidiscip Optim 2011;43:1–16. https://doi.org/10.1007/s00158-010-0594-7.

[35] Bourdin B. Filters in topology optimization. Int J Numer Methods Eng 2001;50:2143–58. https://doi.org/10.1002/nme.116.

[36] Lazarov BS, Sigmund O. Filters in topology optimization based on Helmholtz-type differential equations. Int J Numer Methods Eng 2011;86:765–81. https://doi.org/10.1002/nme.3072.

[37] Wang F, Lazarov BS, Sigmund O. On projection methods, convergence and robust formulations in topology optimization. Struct Multidiscip Optim 2011;43:767–84. https://doi.org/10.1007/s00158-010-0602-y.

[38] Zhang W, Yang W, Zhou J, Li D, Guo X. Structural Topology Optimization Through Explicit Boundary Evolution. J Appl Mech 2016;84. https://doi.org/10.1115/1.4034972.

[39] Zhang W, Yuan J, Zhang J, Guo X. A new topology optimization approach based on Moving Morphable Components (MMC) and the ersatz material model. Struct Multidiscip Optim 2016;53:1243–60. https://doi.org/10.1007/s00158-015-1372-3.

[40] Du Z, Cui T, Liu C, Zhang W, Guo Y, Guo X. An efficient and easy-to-extend Matlab code of the Moving Morphable Component (MMC) method for three-dimensional topology optimization. Struct Multidiscip Optim 2022;65:158. https://doi.org/10.1007/s00158-022-03239-4.

[41] Kreisselmeier G, Steinhauser R. Systematic control design by optimizing a vector performance index. Comput Aided Des Control Syst 1980:113–7. https://doi.org/10.1016/B978-0-08-024488-4.50022-X.

[42] Svanberg K. The method of moving asymptotes—a new method for structural optimization. Int J Numer Methods Eng 1987;24:359–73. https://doi.org/10.1002/nme.1620240207.

[43] Wu J, Clausen A, Sigmund O. Minimum compliance topology optimization of shell–infill composites for additive manufacturing. Comput Methods Appl Mech Eng 2017;326:358–75. https://doi.org/10.1016/j.cma.2017.08.018.

[44] Meng Z, Li G, Wang B, Hao P. A hybrid chaos control approach of the performance measure functions for reliability-based design optimization. Computers &

Structures 2015;146:32-43. https://doi.org/10.1016/j.compstruc.2014.08.011.
[45] Meng Z, Li G, Yang D, Zhan L. A new directional stability transformation method of chaos control for first order reliability analysis. Structural and Multidisciplinary Optimization 2017;55:601-612. https://doi.org/10.1007/s00158-016-1525-z.
[46] Meng Z, Li G, Wang X, Sait S, Yıldız A. A comparative study of metaheuristic algorithms for reliability-based design optimization problems. Archives of Computational Methods in Engineering 2021; 28:1853-1869. https://doi.org/10.1007/s11831-020-09443-z.

**Figures**

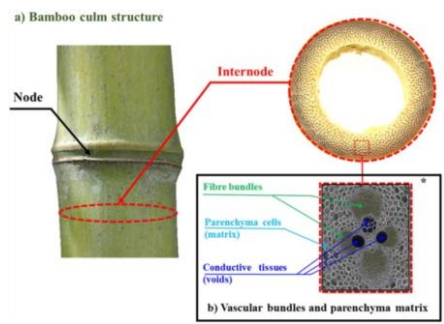 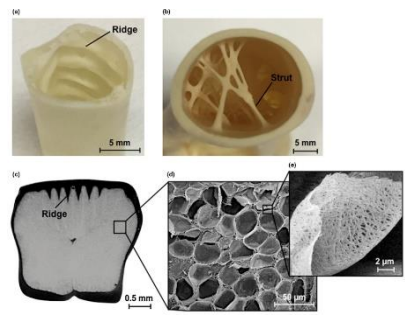

(a) Anatomy of a bamboo [5]   (b) The cross section of a bird's bone [6]

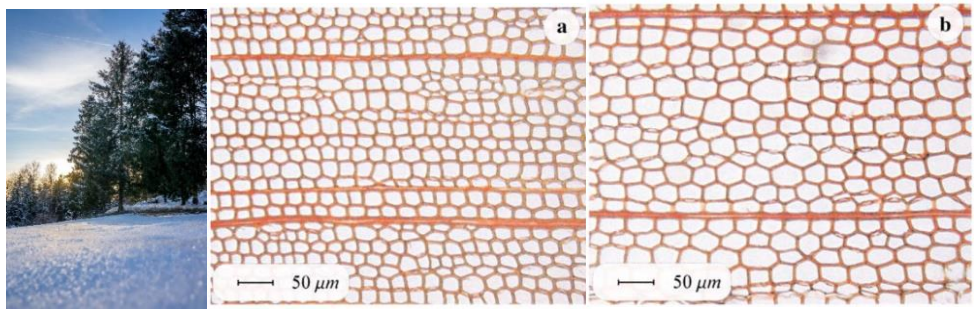

(c) Microscope images of the stem of the Norway spruces [7]

**Fig. 1.** Shell-infill structures in nature

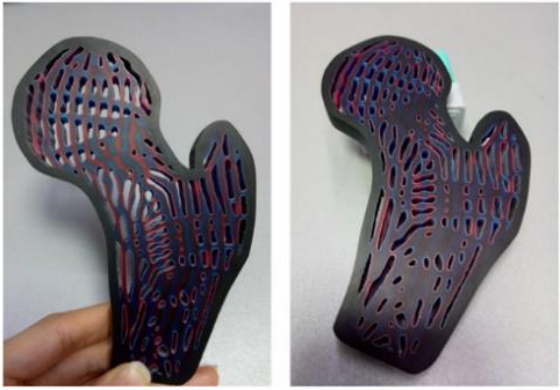 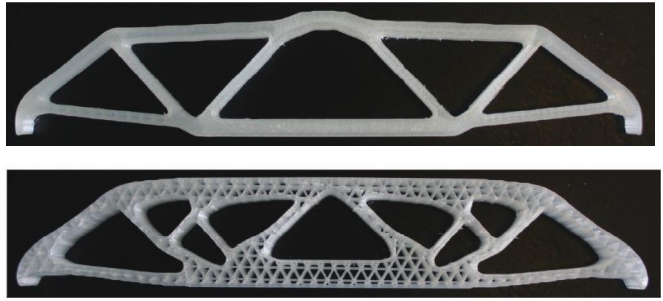

(a) Specimen of the optimized four-phase hip bone [8]     (b) Shell-infill structures fabricated by 3D printers [3]

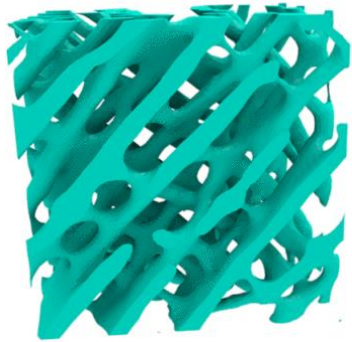 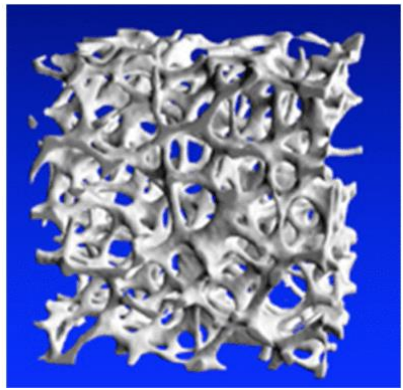

(c) A 3D printed bone-infill structure (left) and a real bone sample from CT scans (right) [4]

**Fig. 2.** The morphology and internal topology of shell-infill composite structures.

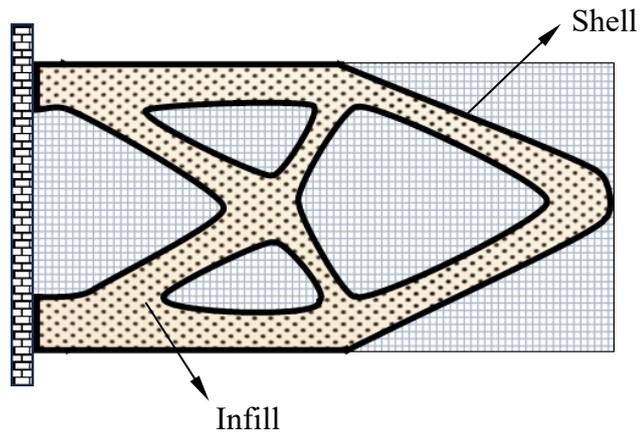

**Fig. 3.** A schematic illustration of the shell-infill structure.

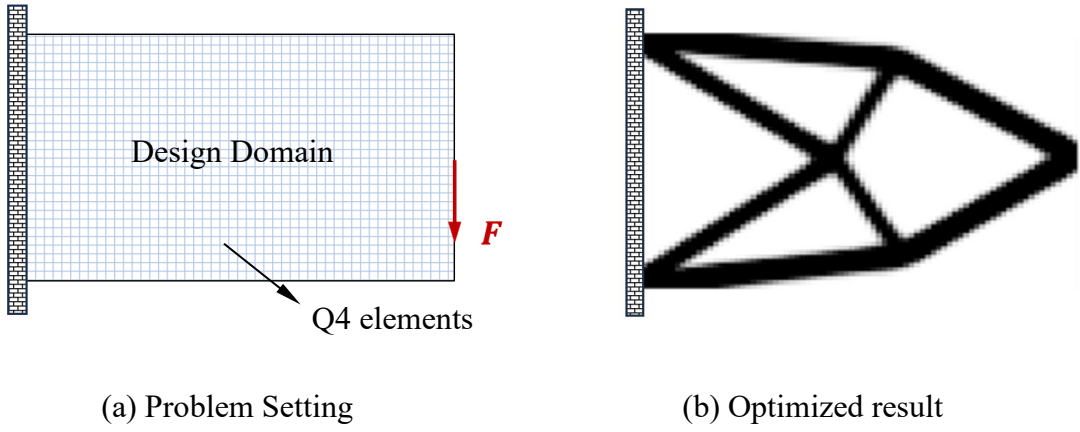

**Fig. 4**. Topology optimization using SIMP method.

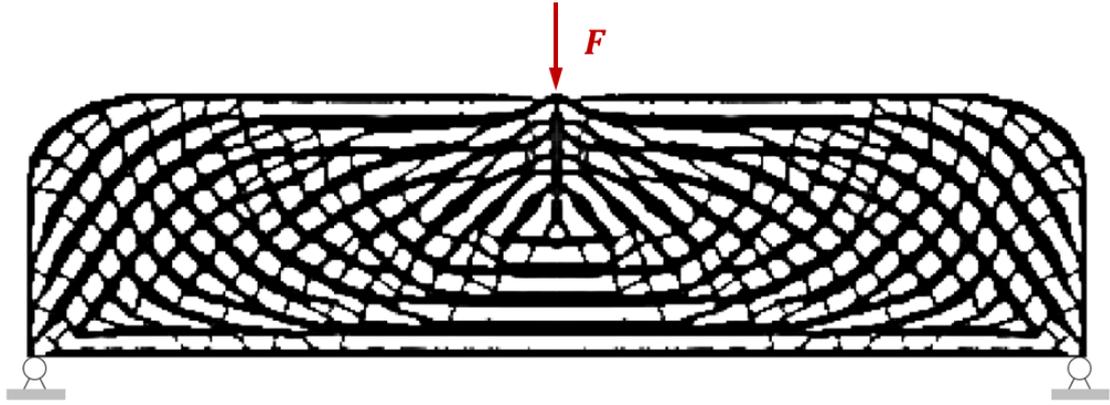

**Fig. 5.** An optimized non-uniform porous infill structure obtained by SIMP method.

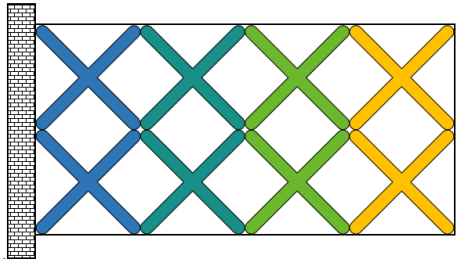
(a) The initial layout of the components

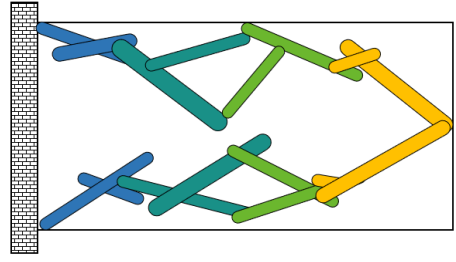
(b) Optimization process

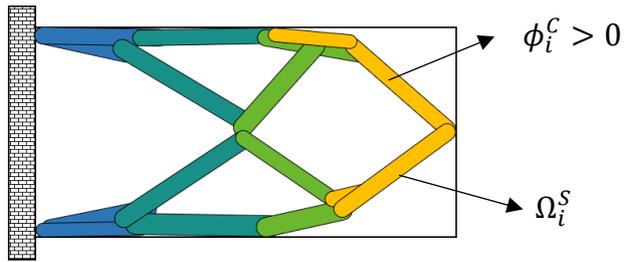
(c) The optimized layout of the components

**Fig. 6.** A schematic illustration of topology optimization based on MMC method.

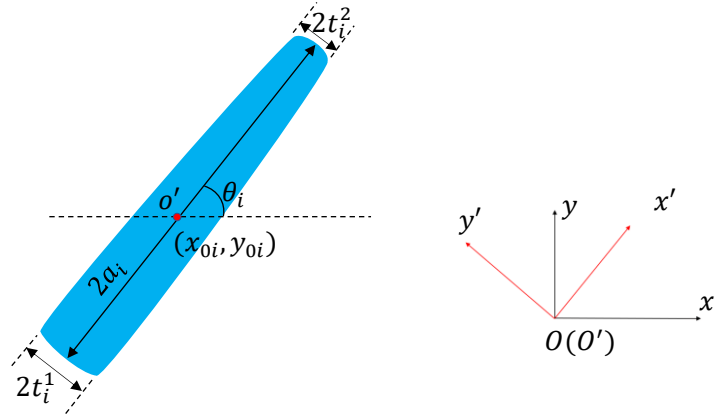

**Fig. 7.** Geometric description of a two-dimensional variable-width component.

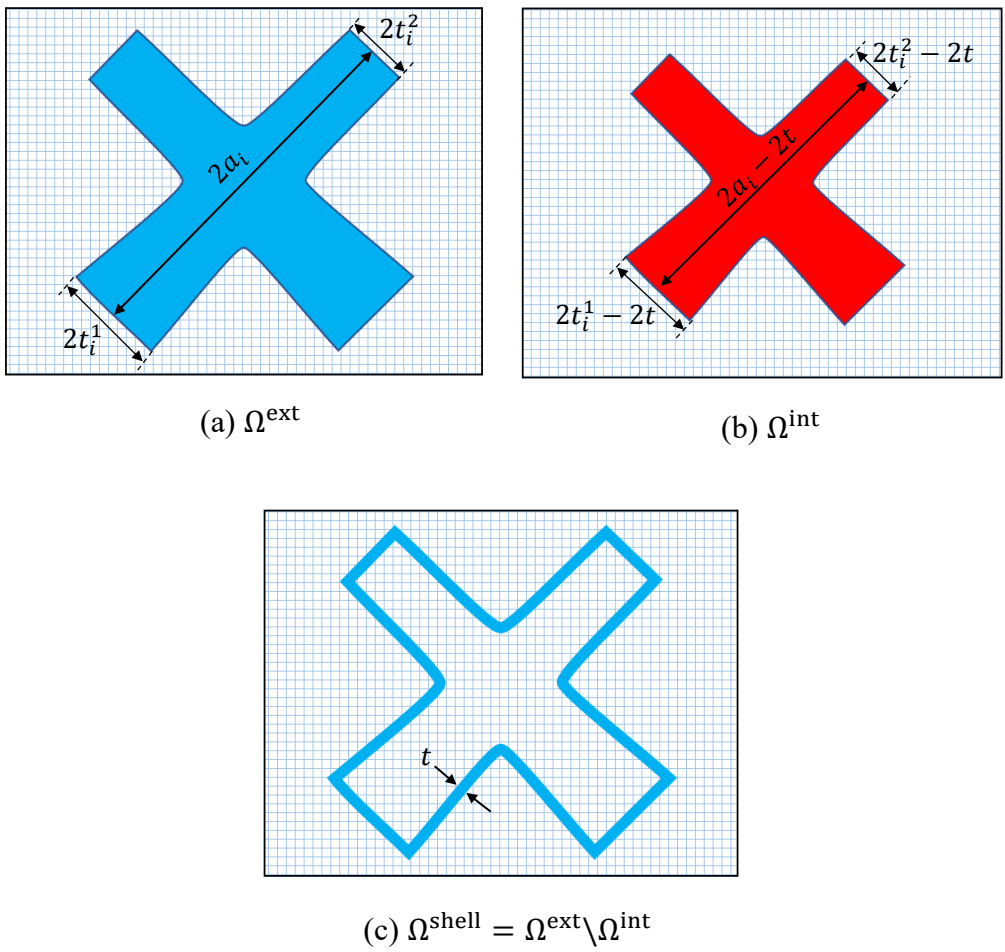

**Fig. 8.** A schematic illustration for the construction of a coating shell.

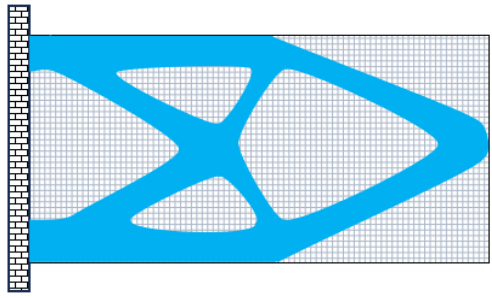 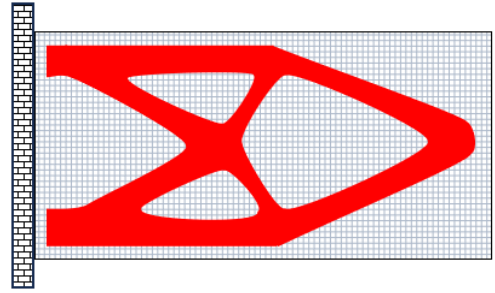

(a) The region $\Omega^{\text{ext}}$ with density vector $\boldsymbol{\rho}^1$

(b) The region $\Omega^{\text{int}}$ with density field $\boldsymbol{\rho}^2$

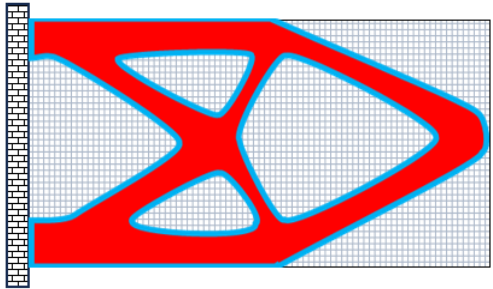 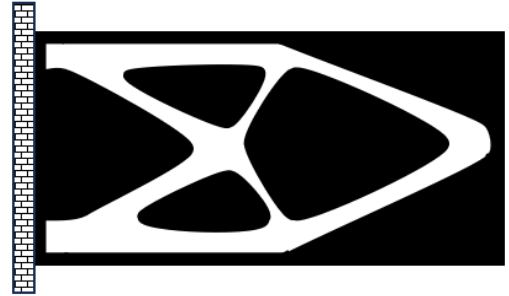

(c) The region $\Omega^{\text{int}} \cup \Omega^{\text{ext}}$

(d) The region $\overline{\Omega^{\text{int}}}$ with density vector $\boldsymbol{I} - \boldsymbol{\rho}^2$

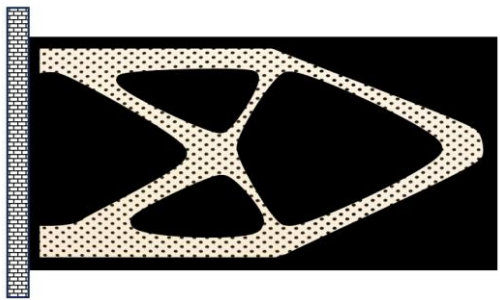 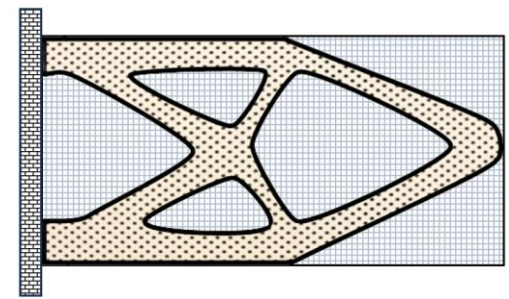

(e) The region $\Omega^1$ with density vector $\boldsymbol{\rho}' = \max(\boldsymbol{\rho}^S, \boldsymbol{I} - \boldsymbol{\rho}^2)$

(f) The region $\Omega^{S-I}$ with the density vector $\boldsymbol{\rho} = \min(\boldsymbol{\rho}', \boldsymbol{\rho}^1)$

**Fig. 9.** A schematic illustration of generating shell-infill structure.

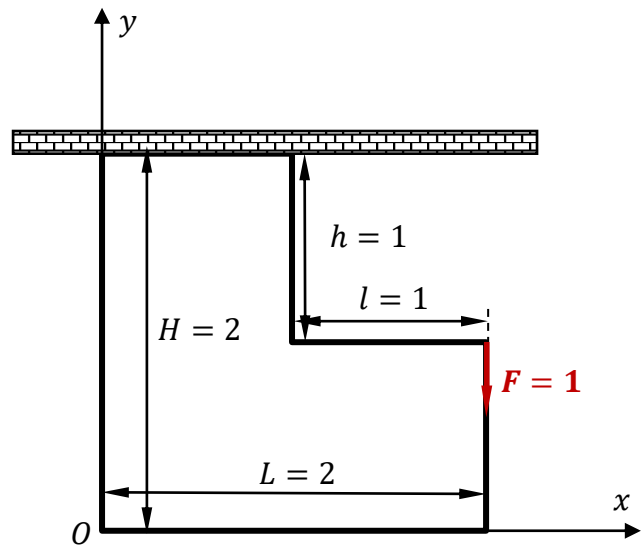

**Fig. 10.** The L-shaped beam example.

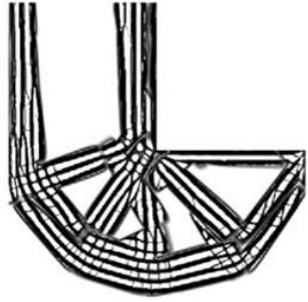

(a) Optimized result without the regularization

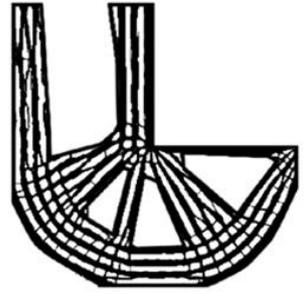

(b) Optimized result with the regularization

**Fig. 11.** Optimization results obtained with and without the regularization term in the objective function.

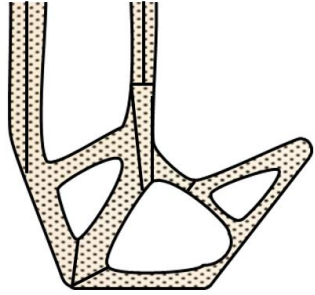

(a) Schematic illustration of the result without the regularization

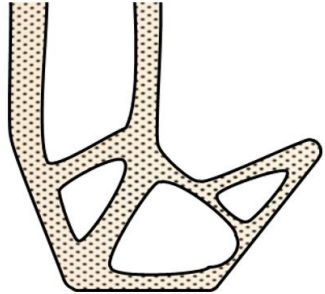

(b) Schematic illustration of the result with the regularization

**Fig. 12.** Schematic illustration of the regularization term promoting the formation of large-area continuous infill structure.

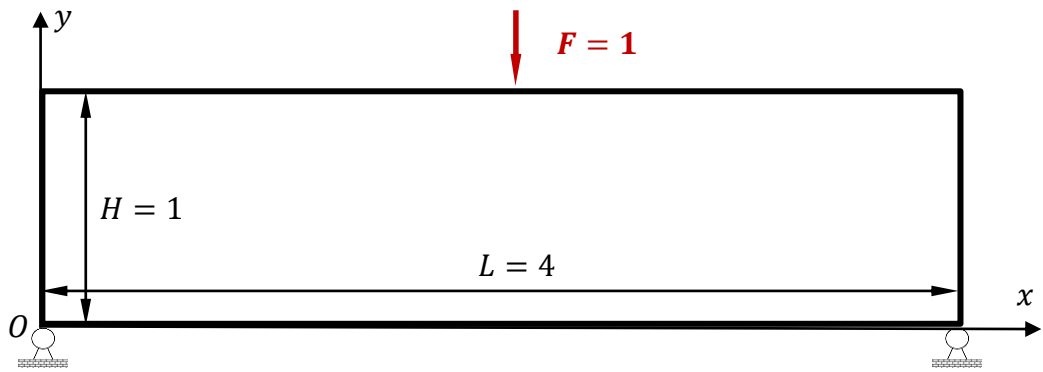

**Fig. 13.** The MBB beam example.

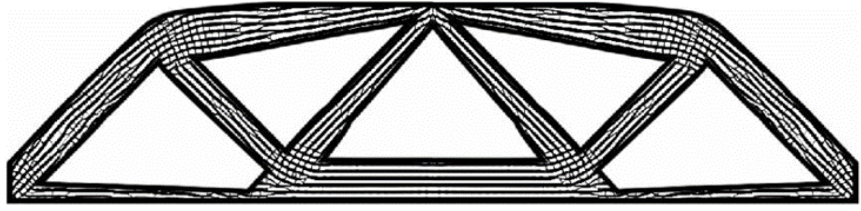

(a) $\bar{V}_c = 0.5$, Compliance = 138.05

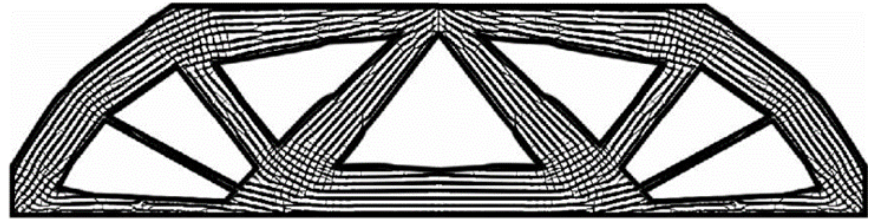

(b) $\bar{V}_c = 0.6$, Compliance = 123.80

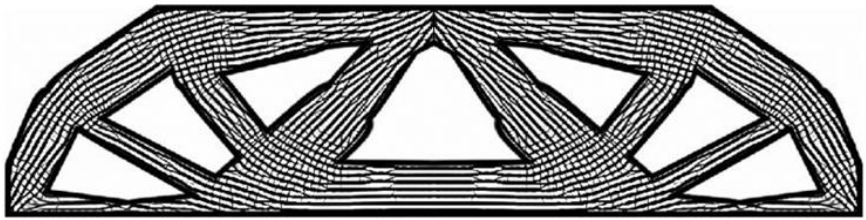

(c) $\bar{V}_c = 0.7$, Compliance = 117.31

**Fig. 14.** The optimized results for the MBB beam example under different component volumes.

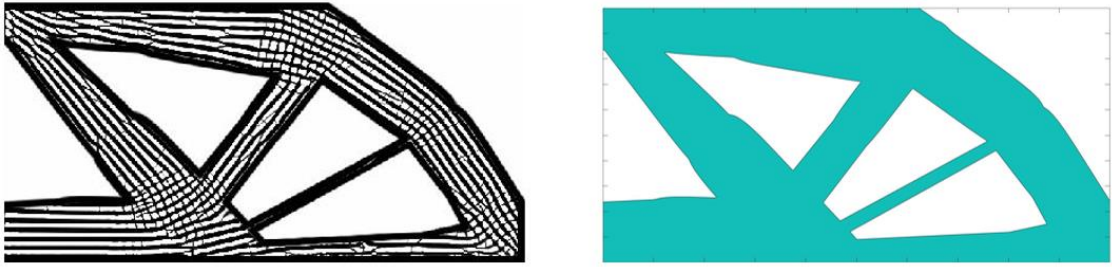

**Fig. 15.** The correlation between the porous infill structure and the base area enclosed by the shell for the MBB beam example.

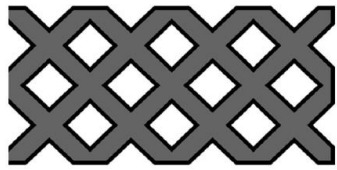 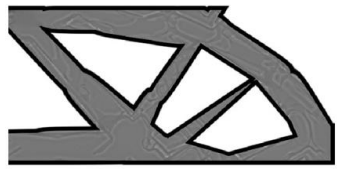 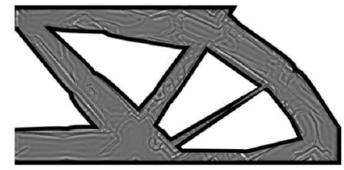

(a) Step 1  (b) Step 200  (c) Step 400

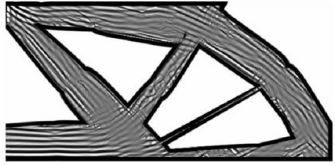 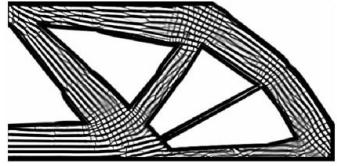 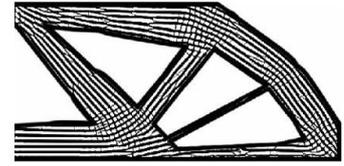

(d) Step 800  (e) Step 1000  (f) Step 1400

**Fig. 16.** The intermediate results of optimization with step 1, 200, 400, 800, 1000, 1400.

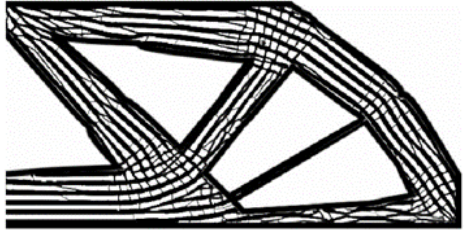
(a) $R_e = 8a$

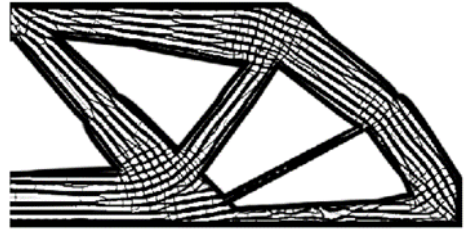
(b) $R_e = 7a$

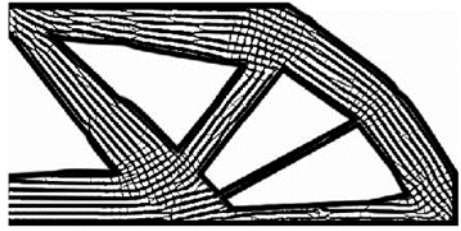
(c) $R_e = 6a$

**Fig. 17.** The optimized results for the MBB beam example under different local volume influence radius.

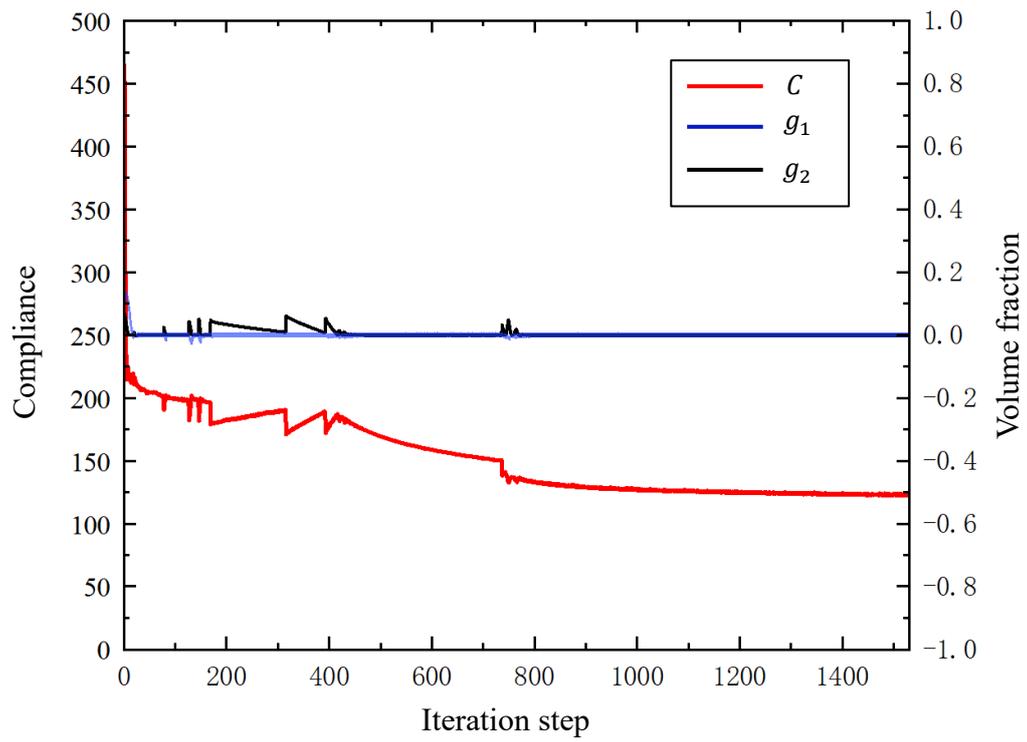

**Fig. 18.** The iteration history of the MBB beam example.

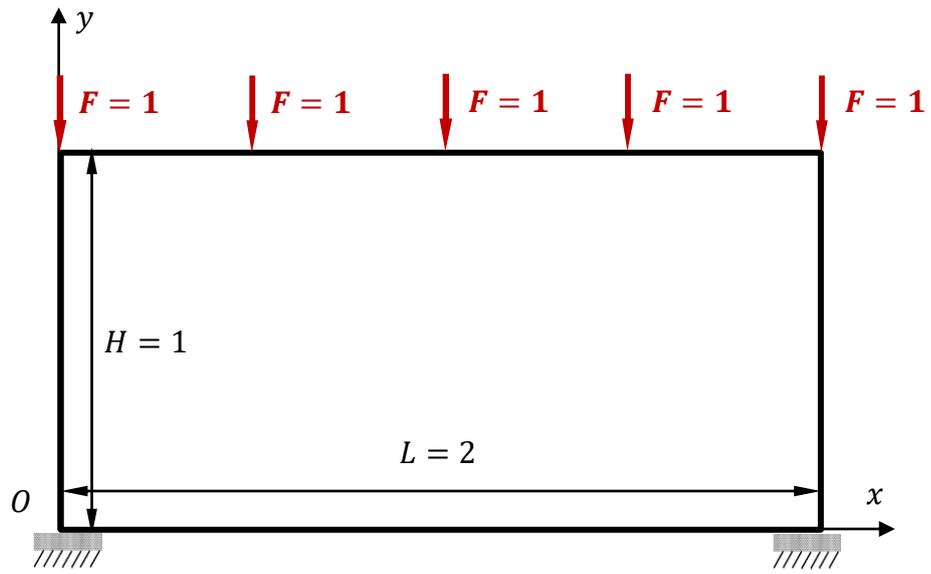

**Fig. 19.** The multi-load beam example.

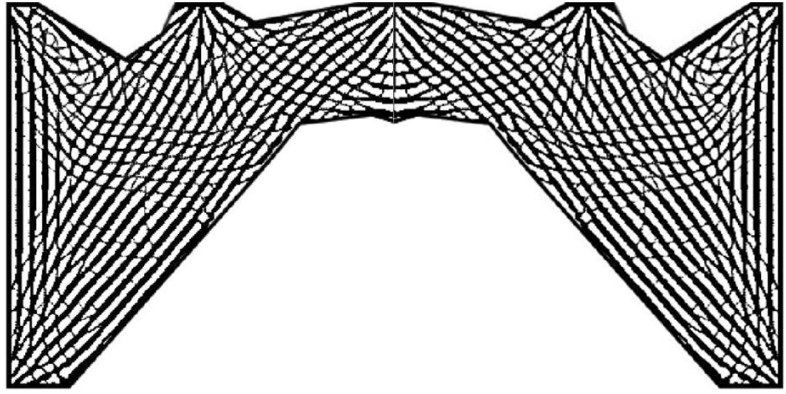

**Fig. 20.** The optimized result for the multi-load example.

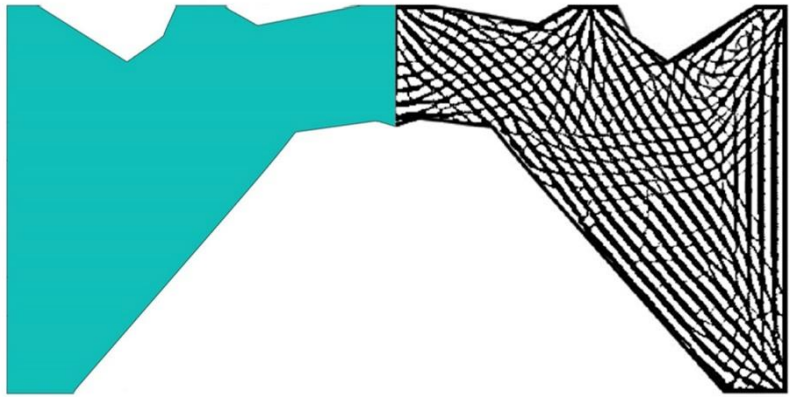

**Fig. 21.** The correlation between the porous infill structure and the base area enclosed by the shell for the multi-load example.